\documentclass[twoside,11pt,reqno]{amsart}

\usepackage{amsmath,amsthm,amssymb,amstext,amsfonts,amscd}
\usepackage{graphicx}
\usepackage{xcolor}
\usepackage{multirow}

\newtheorem{theorem}{Theorem}
\theoremstyle{plain}

\newtheorem{corollary}{Corollary}

\numberwithin{equation}{section}

\begin{document}
\title[\hfil Certain Fractional Kinetic Equations Involving Generalized k-Bessel Function]
{Certain Fractional Kinetic Equations Involving Generalized k-Bessel Function}

\author[ P. Agarwal, S. Jain, A. Atangana,  M. Chand and G. Singh ]{ Praveen
Agarwal, Shilpi Jain, Abdon Atangana, Mehar Chand  and Gurmej Singh}

\bigskip
\address{P. Agarwal:
Department of Mathematics,
Anand International College of Engineering,
Jaipur303012, India}
 \email{goyal.praveen2011@gmail.com}
\address{S. Jain: Department of Mathematics, Poornima College of Engineering, Jaipur-303029, India}
\email{shilpijain1310@gmail.com}

\address{A. Atangana:
Faculty / Fakulteit: Natural and Agricultural Sciences / Natuur-en Landbouwetenskappe
PO Box / Posbus 339, Bloemfontein 9300, Republic of South Africa / Republiek van Suid-Afrika}
 \email{AtanganaA@ufs.ac.za}

\bigskip
\address{M. Chand: Department of Mathematics, Fateh College for Women, Bathinda 151103, India}
\email{mehar.jallandhra@gmail.com}

\bigskip
\address{G. Singh:Department of Mathematics, Mata Sahib Kaur Girls College, Talwandi Sabo, Bathinda-151103 (India)\newline Research Scholar, Department of Mathematics, Singhania University, Pacheri Bari, Jhunjhunu-(India)}
\email{gurmejsandhu11@gmail.com}

\subjclass[2010]{26A33, 33C45, 33C60, 33C70}
\keywords{Gamma function; Beta function; k-Bessel function; Mellin-Barnes type integral; Laguerre polynomials;  Konhauser polynomials}


\thanks{$^{*}$ corresponding author}

\begin{abstract}
We develop a new and further generalized form of the fractional kinetic equation involving generalized k-Bessel function. The manifold generality of the generalized k-Bessel function is discussed in terms of the solution of the fractional
kinetic equation in the present paper. The results obtained here are quite general in nature and capable of yielding a very large number of known and (presumably) new results.
 \end{abstract}

\maketitle

\section{Introduction and Preliminaries}

In recent years, unified integrals involving Special functions attract the attention of the many researchers due to various application point of view(see, \cite{11,12}). In the sequel, Diaz and Pariguan \cite{111} introduced the $k$-Pochhemmer symbol  and $k$-gamma function defined as follows:

\begin{equation}\label{Poch}
\aligned & (\lambda)_{n,k}:
  =\left\{\aligned & \frac{\Gamma_k(\gamma+nk)}{\Gamma_k(\gamma)}  \hskip 41 mm (k\in\mathbb{R};\gamma \in \mathbb{C}\setminus \{0\}) \\
        & \gamma(\gamma+k)...(\gamma+(n-1)k) \hskip20mm (n\in\mathbb{N};\gamma \in \mathbb{C}),
   \endaligned \right.
\endaligned
\end{equation}

They gave the relation with the classical Euler's gamma function(see\cite{18, 10}) as:

\begin{equation}\label{Poch-1}
\aligned \Gamma_k(\gamma)=k^{\frac{\gamma}{k}-1}\Gamma\left(\frac{\gamma}{k}\right)
\endaligned
\end{equation}

Clearly, for $ k=1 $, \eqref{Poch} reduces to the classical  Pochhemmer symbol and Euler's gamma function, respectively (see\cite{116})
\newline

Recently ,Romero et. al.\cite{10} (see, also\cite{118}) introduced the k-Bessel function of the first kind for $\alpha,\lambda,\gamma,\nu\in \mathbb{C} $ and $ \Re(\lambda)>0, \Re(\nu)>0 $ as follows:

\begin{equation}\label{k-bessel}
\aligned J^{(\gamma),(\lambda)}_{k,\nu}(z)=\sum^\infty_{n=0}\frac{(\lambda)_{n,k}}{\Gamma_k(\lambda n+\nu+1)}\frac{(-1)^n\left(\frac{z}{2}\right)^n}{(n!)^2}
\endaligned
\end{equation}

The Fox-Wright function ${}_p\psi_q(z)$ with $p$ numerator and $q$ denominators, such that $a_i,b_j\in\mathbb{C} (i=1,...,p;j=1,...,q)$ is defined by (see, for detail\cite{112}):

\begin{eqnarray}\label{wright-function}{}_p\psi_q(z)={}_p\psi_q\left[
\left. \begin{array}{cc} (a_i,\alpha_i)_{1,p}
\\(b_j,\beta_j)_{1,q}	
\end{array}\right|z
\right]=\sum^\infty_{n=0}\frac{\prod^p_{i=1}\Gamma(a_i+\alpha_in)}{\prod^q_{j=1}\Gamma(b_j+\beta_Jn)}\frac{z^n}{n!} \end{eqnarray}

under the condition

\begin{eqnarray} \sum^{q}_{j=1}\beta_j-\sum^p_{i=1}\alpha_i>-1 \end{eqnarray}

In particular, when $a_i=b_j=1(i=1,...,p;j=1,...,q),$ immediate reduces to the generalized hypergeometric function ${}_pF_q (p,q \in \mathbb{N}_0) $  (see, for details\cite{15}):

\begin{eqnarray}\label{wright-function-case-1}{}_p\psi_q(z)={}_p\psi_q\left[
\left. \begin{array}{cc} (a_i,1)_{1,p}
\\(b_j,1)_{1,q}	
\end{array}\right|z
\right]=\frac{\prod^p_{i=1}\Gamma(a_i)}{\prod^q_{j=1}\Gamma(b_j)}{}_pF_q\left[\begin{array}{cc} a_1,...,a_p;
\\b_1,...,b_q;	
\end{array} z\right] \end{eqnarray}


In terms of the $k$-Pochhamer symbol $(\gamma)_{n,k}$ defined by \eqref{Poch}, we introduce more generalized form of k-Bessel function $\omega^{\gamma,\lambda}_{k,\nu,b,c}(z)$ as follows:

\begin{eqnarray}\label{k-bessel-function}\omega^{\gamma,\lambda}_{k,\nu,b,c}(z)=\sum^\infty_{n=0}\frac{(-1)^nc^n(\gamma)_{n,k}}{\Gamma_k(\nu+\lambda n+\frac{b+1}{2})}\frac{\left(\frac{z}{2}\right)^{\nu+2n}}{(n!)^2} \end{eqnarray}
where $\alpha,\lambda,\gamma,\nu,c,b\in\mathbb{C}$ and $\Re(\lambda)>0,\Re(\nu)>0$.
\newline

\noindent The importance of fractional differential equations in the field of applied science has gained more attention not only in mathematics but also in physics, dynamical systems, control systems and engineering, to create the mathematical model of many physical phenomena.  Especially, the kinetic equations describe the continuity of motion of substance. The extension and generalization of fractional kinetic equations involving many fractional operators were found \cite{5,6,7,8,9,10,11,12,13,14,15,16,17,18}.
\newline

In view of the effectiveness and a great importance of the kinetic equation in certain astrophysical problems the authors develop a further generalized form of the fractional kinetic equation involving generalized k-Bessel function.
\newline

The fractional differential equation between rate of change of the reaction was established by Haubold and Mathai\cite{7}, the destruction rate and the production rate are given as follows:

\begin{eqnarray}\label{D-1} \frac{dN}{dt}=-d(N_{t})+p(N_{t}), \end{eqnarray}

    where $ N=N(t) $ the rate of reaction,  $ d=d(N) $ the rate of destruction,   $ p=p(N) $ the rate of production and $N_{t} $ denotes the function defined by   $N_{t}(t^{*})= N(t-t^{*}), t^{*}>0 $
\newline

   The special case of \eqref{D-1} for spatial fluctuations and inhomogeneities in $ N(t) $  the quantities are neglected , that is the equation
\begin{eqnarray}\aligned &\label{D-2} \frac{dN}{dt}=-c_{i}N_{i}(t),\endaligned \end{eqnarray}

with the initial condition that $ N_{i}(t=0)=N_{0} $ is the number density of the species $ i $ at time $ t=0 $ and $ c_{i}>0 $. If we remove the index $ i $ and integrate the standard kinetic equation \eqref{D-2}, we have

\begin{eqnarray}\aligned &\label{D-3} N(t)-N_{0}=-c{}_{0}D^{-1}_{t}N(t)\endaligned \end{eqnarray}

where $ {}_{0}D^{-1}_{t} $ is the special case of the Riemann-Liouville integral operator $ {}_{0}D^{-\nu}_{t} $ defined as

\begin{eqnarray}\aligned &\label{D-4} {}_{0}D^{-\nu}_{t}f(t)=\frac{1}{\Gamma(\nu)}\int_{0}^{t}\left(t-s\right)^{\nu-1}f(s)ds,\hskip 6mm  (t>0,R(\nu)>0)\endaligned \end{eqnarray}

     The fractional generalization of the standard kinetic equation\eqref{D-3} is given by Haubold and Mathai\cite{7} as follows:
		
\begin{eqnarray}\label{D-5} N(t)-N_{0}=-c^{\nu}{}_{0}D^{-1}_{t}N(t) \end{eqnarray}

 and obtained the solution of \eqref{D-5} as follows:

\begin{eqnarray}\aligned &\label{D-6} N(t)=N_{0}\sum _{k=0}^{\infty }\frac{(-1)^{k}}{\Gamma\left(\nu k+1\right)}\left(ct\right)^{\nu k}\endaligned \end{eqnarray}

Further, \cite{11} considered the the following fractional kinetic equation:

\begin{eqnarray}\aligned &\label{D-7} N(t)-N_{0}f(t)=-c^{\nu}{}_{0}D^{-\nu}_{t}N(t),\hskip 6mm (\Re(v)>0),\endaligned\end{eqnarray}

where $ N(t) $ denotes the number density of a given species at time $ t $, $ N_{0}=N(0) $ is the number density of that species at time $t=0 $, $ c $ is a constant and $ f \in \mathcal{L}(0,\infty) $.
\newline

 By applying the Laplace transform to \eqref{D-7} (see\cite{17}),
\begin{eqnarray}\aligned &\label{D-8} L\left\{ N(t);p\right\}=N_{0}\frac{F(p)}{1+c^{\nu}p^{-\nu}}=N_{0}\left(\sum _{n=0}^{\infty }(-c^{\nu})^{n}p^{-\nu n}\right)F(p),\\&\hskip 6mm \left(n\in N_{0},\left|\frac{c}{p}\right|<1\right)\endaligned \end{eqnarray}

where the Laplace transform \cite{19} is given by

\begin{eqnarray}\aligned &\label{D-9} F(p)=L\left\{ N(t);p\right\} =\int_{0}^{\infty} e^{-pt}f(t)dt,\hskip 6mm (\mathcal{R}(p)>0).\endaligned  \end{eqnarray}

The objective of this paper is to derive the solution of the fractional kinetic equation involving generalized k-Bessel function. The results obtained in terms of Mittag-Leffler function are rather  general in nature and can easily construct various known and new fractional kinetic equations.

 \section{Solution of generalized fractional kinetic equations}
In this section, we will investigate the solution of the generalized fractional kinetic equations by considering generalized k-Bessel function. The results are as follows.

\begin{theorem}\label{Th-1} If $ d>0, \nu>0, \lambda,\gamma,\mu,c,b\in\mathbb{C}, k\in\mathbb{N}$ and $\Re(\lambda)>0,\Re(\mu)>0$ then the solution of the equation
\begin{eqnarray}\aligned &\label{Th-1-1} N(t)-N_{0}\omega^{\gamma,\lambda}_{k,\mu,b,c}(t)=-d^{\nu}{}_{0}D^{-\nu}_{t}N(t)\endaligned\end{eqnarray}
is given by the following formula
\begin{eqnarray}\aligned &\label{Th-1-2} N(t)=N_{0}\sum^\infty_{n=0}\frac{(-1)^nc^n(\gamma)_{n,k}}{\Gamma_k(\mu+\lambda n+\frac{b+1}{2})}\frac{\Gamma(\mu+2n+1)}{(n!)^2}\left(\frac{t}{2}\right)^{\mu+2n} E_{\nu,\hskip .5mm\mu+2n+1}(-d^{\nu}t^{\nu}),\endaligned \end{eqnarray}
where the generalized Mittag-Leffler function $ E_{\alpha,\beta}(x) $ is given by \cite{20}

\begin{eqnarray}\label{2.3}\aligned &\label{ML} E_{\alpha,\beta}(x)=\sum _{n=0}^{\infty }\frac{(x)^{n}}{\Gamma\left(\alpha n+\beta\right)}.\endaligned \end{eqnarray}

\end{theorem}

Proof: the Laplace transform of Riemann-Liouville fractional integral operator is given by \cite{21} \cite{22}

\begin{eqnarray}\aligned &\label{Th-1-3 } L\left\{{}_{0}D^{-\nu}_{t}f(t);p\right\}=p^{-\nu}F(p)\endaligned \end{eqnarray}

where $ F(p) $ is defined in \eqref{D-9}.Now ,applying the Laplace transform to the both sides of \eqref{Th-1-1} gives

\begin{eqnarray}\aligned &\label{Th-1-4}L\left\{ N(t);p\right\}=N_{0}L\left\{ \omega^{\gamma,\lambda}_{k,\mu,b,c}(t);p\right\}-d^{\nu}L\left\{{}_{0}D^{-\nu}_{t}N(t);p\right\}\endaligned \end{eqnarray}

\begin{eqnarray}\aligned &\label{Th-1-5} N(p)=N_{0}\left(\int^{\infty}_{0}e^{-pt}\sum^\infty_{n=0}\frac{(-1)^nc^n(\gamma)_{n,k}}{\Gamma_k(\mu+\lambda n+\frac{b+1}{2})}\frac{1}{(n!)^2}\left(\frac{t}{2}\right)^{\mu+2n}dt\right)\\&\hskip 6mm-d^{\nu}p^{-\nu}N(p)\endaligned\end{eqnarray}

\begin{eqnarray}\aligned &\label{Th-1-6} N(p)+d^{\nu}p^{-\nu}N(p)=N_{0}\sum^\infty_{n=0}\frac{(-1)^nc^n(\gamma)_{n,k}}{\Gamma_k(\mu+\lambda n+\frac{b+1}{2})}\frac{1}{(n!)^2}\left(\frac{1}{2}\right)^{\mu+2n}\\&\hskip 6mm\times\int^{\infty}_{0}e^{-pt}t^{\mu+2n}dt\endaligned \end{eqnarray}

\begin{eqnarray}\aligned &\label{Th-1-7} =N_{0}\sum^\infty_{n=0}\frac{(-1)^nc^n(\gamma)_{n,k}}{\Gamma_k(\mu+\lambda n+\frac{b+1}{2})}\frac{1}{(n!)^2}\left(\frac{1}{2}\right)^{\mu+2n}\\&\hskip 6mm\times\frac{\Gamma(\mu+2n+1)}{p^{\mu+2n+1}}\endaligned\end{eqnarray}

\begin{eqnarray}\aligned &\label{Th-1-8} N(p)=N_{0}\sum^\infty_{n=0}\frac{(-1)^nc^n(\gamma)_{n,k}}{\Gamma_k(\mu+\lambda n+\frac{b+1}{2})}\frac{\Gamma(\mu+2n+1)}{(n!)^2}\left(\frac{1}{2}\right)^{\mu+2n}\\&\hskip 6mm\times\left\{p^{-(\mu+2n+1)}\sum _{r=0}^{\infty }\left[-\left(\frac{p}{d}\right)^{-\nu}\right]^{r}\right\}\endaligned \end{eqnarray}

Taking Laplace inverse of \eqref{Th-1-8},and by using

\begin{eqnarray}\aligned &\label{Th-1-9} L^{-1}\left\{p^{-\nu};t\right\}=\frac{t^{\nu-1}}{\Gamma(\nu)},(R(\nu)>0)\endaligned\end{eqnarray}

we have
\begin{eqnarray}\aligned &\label{Th-1-10} L^{-1}\left\{N(p)\right\}=N_{0}\sum^\infty_{n=0}\frac{(-1)^nc^n(\gamma)_{n,k}}{\Gamma_k(\mu+\lambda n+\frac{b+1}{2})}\frac{\Gamma(\mu+2n+1)}{(n!)^2}\left(\frac{1}{2}\right)^{\mu+2n}\\&\hskip 6mm\times L^{-1}\left\{\sum _{r=0}^{\infty }(-1)^{r}d^{\nu r}p^{-(\mu+2n+1+\nu r)}\right\}\endaligned \end{eqnarray}

\begin{eqnarray}\aligned &\label{Th-1-11} N(t)=N_{0}\sum^\infty_{n=0}\frac{(-1)^nc^n(\gamma)_{n,k}}{\Gamma_k(\mu+\lambda n+\frac{b+1}{2})}\frac{\Gamma(\mu+2n+1)}{(n!)^2}\left(\frac{1}{2}\right)^{\mu+2n}\\&\hskip 6mm\times\left\{\sum _{r=0}^{\infty }(-1)^{r}d^{\nu r}\frac{t^{(\mu+2n+\nu r)}}{\Gamma\left(\mu+2n+\nu r+1\right)}\right\}\endaligned\end{eqnarray}

\begin{eqnarray}\aligned &\label{Th-1-12} =N_{0}\sum^\infty_{n=0}\frac{(-1)^nc^n(\gamma)_{n,k}}{\Gamma_k(\mu+\lambda n+\frac{b+1}{2})}\frac{\Gamma(\mu+2n+1)}{(n!)^2}\left(\frac{t}{2}\right)^{\mu+2n}\\&\hskip 6mm\times\left\{\sum _{r=0}^{\infty }(-1)^{r}d^{\nu r}\frac{t^{(\nu r)}}{\Gamma\left(\mu+2n+\nu r+1\right)}\right\}\endaligned \end{eqnarray}

\begin{eqnarray}\aligned &\label{Th-1-13} N(t)=N_{0}\sum^\infty_{n=0}\frac{(-1)^nc^n(\gamma)_{n,k}}{\Gamma_k(\mu+\lambda n+\frac{b+1}{2})}\frac{\Gamma(\mu+2n+1)}{(n!)^2}\left(\frac{t}{2}\right)^{\mu+2n} E_{\nu,\mu+2n+1}(-d^{\nu}t^{\nu})\endaligned \end{eqnarray}

\begin{theorem}\label{Th-2} If $ d>0, \nu>0, \lambda,\gamma,\mu,c,b\in\mathbb{C}, k\in\mathbb{N}$ and $\Re(\lambda)>0,\Re(\mu)>0$ then the solution of the equation

\begin{eqnarray}\aligned &\label{Th-2-1} N(t)=N_{0}\omega^{\gamma,\lambda}_{k,\mu,b,c}(d^{\nu}t^{\nu})-d^{\nu}{}_{0}D^{-\nu}_{t}N(t)\endaligned\end{eqnarray}

is given by the following formula

\begin{eqnarray}\aligned &\label{Th-2-2} N(t)=N_{0}\sum^\infty_{n=0}\frac{(-1)^nc^n(\gamma)_{n,k}}{\Gamma_k(\mu+\lambda n+\frac{b+1}{2})}\frac{\Gamma(\nu(\mu+2n)+1)}{(n!)^2}\left(\frac{d^{\nu}t^{\nu}}{2}\right)^{\mu+2n}\\&\hskip 10mm\times E_{\nu,\hskip .5mm\nu(\mu+2n)+1}(-d^{\nu}t^{\nu}),\endaligned\end{eqnarray}
where $ E_{\nu,\nu(\mu+2n)+1}(.) $ is the generalized Mittag-Leffler function defined in equation \eqref{2.3}.
\end{theorem}

\begin{theorem}\label{Th-3} If $ a>0, d>0, \nu>0; a\neq d; \lambda,\gamma,\mu,c,b\in\mathbb{C}, k\in\mathbb{N}$ and $\Re(\lambda)>0,\Re(\mu)>0$ then the solution of the equation

\begin{eqnarray}\aligned &\label{Th-3-1} N(t)=N_{0}\omega^{\gamma,\lambda}_{k,\mu,b,c}(d^{\nu}t^{\nu})-a^{\nu}{}_{0}D^{-\nu}_{t}N(t)\endaligned\end{eqnarray}

is given by the following formula

\begin{eqnarray}\aligned &\label{Th-3-2} N(t)=N_{0}\sum^\infty_{n=0}\frac{(-1)^nc^n(\gamma)_{n,k}}{\Gamma_k(\mu+\lambda n+\frac{b+1}{2})}\frac{\Gamma(\nu(\mu+2n)+1)}{(n!)^2}\left(\frac{d^{\nu}t^{\nu}}{2}\right)^{\mu+2n}\\&\hskip 10mm\times E_{\nu,\hskip .5mm\nu(\mu+2n)+1}(-a^{\nu}t^{\nu}),\endaligned\end{eqnarray}
where $ E_{\nu,\nu(\mu+2n)+1}(.) $ is the generalized Mittag-Leffler function defined in equation \eqref{2.3}.
\end{theorem}

Proof: The proof of theorem 2 and 3 would run parallel to those of theorem 1.

\section{Special Cases}
If we choose $b=c=1$ then generalized k-Bessel function reduced to the following form:

\begin{eqnarray}\omega^{\gamma,\lambda}_{k,\mu,1,1}(z)=\left(\frac{z}{2}\right)^{\mu}\sum^\infty_{n=0}\frac{(-1)^n(\gamma)_{n,k}}{\Gamma_k(\lambda n+\mu+1)}\frac{\left(\frac{z^2}{4}\right)^{n}}{(n!)^2}=\left(\frac{z}{2}\right)^{\mu}J^{(\gamma),(\lambda)}_{k,\mu}\left(\frac{z^2}{2}\right) ,\end{eqnarray}

where $\lambda,\gamma,\mu,\in\mathbb{C}$ and $\Re(\lambda)>0,\Re(\mu)>0$.
\newline

Then the Theorems 1, 2 and 3 reduced to the following the form:

\begin{corollary} If $ d>0, \nu>0, \lambda,\gamma,\mu\in\mathbb{C}, k\in\mathbb{N}$ and $\Re(\lambda)>0,\Re(\mu)>0$ then the solution of the equation
\begin{eqnarray}\aligned & N(t)-N_{0}\left(\frac{t}{2}\right)^{\mu}J^{(\gamma),(\lambda)}_{k,\mu}\left(\frac{t^2}{2}\right)=-d^{\nu}{}_{0}D^{-\nu}_{t}N(t)\endaligned\end{eqnarray}
is given by the following formula
\begin{eqnarray}\aligned &  N(t)=N_{0}\sum^\infty_{n=0}\frac{(-1)^n(\gamma)_{n,k}}{\Gamma_k(\mu+\lambda n+1)}\frac{\Gamma(\mu+2n+1)}{(n!)^2}\left(\frac{t}{2}\right)^{\mu+2n} E_{\nu,\hskip .5mm\mu+2n+1}(-d^{\nu}t^{\nu}),\endaligned \end{eqnarray}
where $ E_{\nu,\mu+2n+1}(.) $ is the generalized Mittag-Leffler function defined in equation \eqref{2.3}.
\end{corollary}

\begin{corollary} If $ d>0, \nu>0, \lambda,\gamma,\mu\in\mathbb{C}, k\in\mathbb{N}$ and $\Re(\lambda)>0,\Re(\mu)>0$ then the solution of the equation

\begin{eqnarray}\aligned & N(t)=N_{0}\left(\frac{d^{\nu}t^{\nu}}{2}\right)^{\mu}J^{(\gamma),(\lambda)}_{k,\mu}\left(\frac{(d^{\nu}t^{\nu})^2}{2}\right)-d^{\nu}{}_{0}D^{-\nu}_{t}N(t)\endaligned\end{eqnarray}

is given by the following formula

\begin{eqnarray}\aligned & N(t)=N_{0}\sum^\infty_{n=0}\frac{(-1)^n(\gamma)_{n,k}}{\Gamma_k(\mu+\lambda n+1)}\frac{\Gamma(\nu(\mu+2n)+1)}{(n!)^2}\left(\frac{d^{\nu}t^{\nu}}{2}\right)^{\mu+2n}\\&\hskip 10mm\times E_{\nu,\hskip .5mm\nu(\mu+2n)+1}(-d^{\nu}t^{\nu}),\endaligned\end{eqnarray}
where $ E_{\nu,\nu(\mu+2n)+1}(.) $ is the generalized Mittag-Leffler function defined in equation \eqref{2.3}.
\end{corollary}

\begin{corollary} If $ a>0, d>0, \nu>0;a\neq d; \lambda,\gamma,\mu\in\mathbb{C}, k\in\mathbb{N}$ and $\Re(\lambda)>0,\Re(\mu)>0$ then the solution of the equation

\begin{eqnarray}\aligned &  N(t)=N_{0}\left(\frac{d^{\nu}t^{\nu}}{2}\right)^{\mu}J^{(\gamma),(\lambda)}_{k,\mu}\left(\frac{(d^{\nu}t^{\nu})^2}{2}\right)-a^{\nu}{}_{0}D^{-\nu}_{t}N(t)\endaligned\end{eqnarray}

is given by the following formula

\begin{eqnarray}\aligned &  N(t)=N_{0}\sum^\infty_{n=0}\frac{(-1)^n(\gamma)_{n,k}}{\Gamma_k(\mu+\lambda n+1)}\frac{\Gamma(\nu(\mu+2n)+1)}{(n!)^2}\left(\frac{d^{\nu}t^{\nu}}{2}\right)^{\mu+2n}\\&\hskip 10mm\times E_{\nu,\hskip .5mm\nu(\mu+2n)+1}(-a^{\nu}t^{\nu}),\endaligned\end{eqnarray}
where $ E_{\nu,\nu(\mu+2n)+1}(.) $ is the generalized Mittag-Leffler function defined in equation \eqref{2.3}.
\end{corollary}


If we choose $b=-1,c=1$ then generalized k-Bessel function reduced to the k-Wright function \cite{23} associated with the following relation:

\begin{eqnarray}\omega^{\gamma,\lambda}_{k,\mu,-1,1}(z)=\left(\frac{z}{2}\right)^{\mu}\sum^\infty_{n=0}\frac{(-1)^n(\gamma)_{n,k}}{\Gamma_k(\lambda n+\mu)}\frac{\left(\frac{z^2}{4}\right)^{n}}{(n!)^2}=\left(\frac{z}{2}\right)^{\mu}W^{\gamma}_{k,\lambda,\mu}\left(\frac{-z^2}{2}\right) \end{eqnarray}

where $\lambda,\gamma,\mu,\in\mathbb{C}$ and $\Re(\lambda)>0,\Re(\mu)>0$.
\newline

Then the Theorems 1, 2 and 3 reduced to the following the form:

\begin{corollary} If $ d>0, \nu>0, \lambda,\gamma,\mu\in\mathbb{C}, k\in\mathbb{N}$ and $\Re(\lambda)>0,\Re(\mu)>0$ then the solution of the equation
\begin{eqnarray}\aligned & N(t)-N_{0}\left(\frac{t}{2}\right)^{\mu}W^{\gamma}_{k,\lambda,\mu}\left(\frac{-t^2}{2}\right)=-d^{\nu}{}_{0}D^{-\nu}_{t}N(t)\endaligned\end{eqnarray}
is given by the following formula
\begin{eqnarray}\aligned &  N(t)=N_{0}\sum^\infty_{n=0}\frac{(-1)^n(\gamma)_{n,k}}{\Gamma_k(\mu+\lambda n+1)}\frac{\Gamma(\mu+2n+1)}{(n!)^2}\left(\frac{t}{2}\right)^{\mu+2n} E_{\nu,\hskip .5mm\mu+2n+1}(-d^{\nu}t^{\nu}),\endaligned \end{eqnarray}
where $ E_{\nu,\mu+2n+1}(.) $ is the generalized Mittag-Leffler function defined in equation \eqref{2.3}.
\end{corollary}

\begin{corollary} If $ d>0, \nu>0, \lambda,\gamma,\mu\in\mathbb{C}, k\in\mathbb{N}$ and $\Re(\lambda)>0,\Re(\mu)>0$ then the solution of the equation

\begin{eqnarray}\aligned & N(t)=N_{0}\left(\frac{d^{\nu}t^{\nu}}{2}\right)^{\mu}W^{\gamma}_{k,\lambda,\mu}\left(\frac{-(d^{\nu}t^{\nu})^2}{2}\right)-d^{\nu}{}_{0}D^{-\nu}_{t}N(t)\endaligned\end{eqnarray}

is given by the following formula

\begin{eqnarray}\aligned & N(t)=N_{0}\sum^\infty_{n=0}\frac{(-1)^n(\gamma)_{n,k}}{\Gamma_k(\mu+\lambda n+1)}\frac{\Gamma(\nu(\mu+2n)+1)}{(n!)^2}\left(\frac{d^{\nu}t^{\nu}}{2}\right)^{\mu+2n}\\&\hskip 10mm\times E_{\nu,\hskip .5mm\nu(\mu+2n)+1}(-d^{\nu}t^{\nu}),\endaligned\end{eqnarray}
where $ E_{\nu,\nu(\mu+2n)+1}(.) $ is the generalized Mittag-Leffler function defined in equation \eqref{2.3}.

\end{corollary}

\begin{corollary} If $ a>0, d>0, \nu>0 ;a\neq d; \lambda,\gamma,\mu\in\mathbb{C}, k\in\mathbb{N}$ and $\Re(\lambda)>0,\Re(\mu)>0$ then the solution of the equation

\begin{eqnarray}\aligned &  N(t)=N_{0}\left(\frac{d^{\nu}t^{\nu}}{2}\right)^{\mu}W^{\gamma}_{k,\lambda,\mu}\left(-\frac{(d^{\nu}t^{\nu})^2}{2}\right)-a^{\nu}{}_{0}D^{-\nu}_{t}N(t)\endaligned\end{eqnarray}

is given by the following formula

\begin{eqnarray}\aligned &  N(t)=N_{0}\sum^\infty_{n=0}\frac{(-1)^n(\gamma)_{n,k}}{\Gamma_k(\mu+\lambda n+1)}\frac{\Gamma(\nu(\mu+2n)+1)}{(n!)^2}\left(\frac{d^{\nu}t^{\nu}}{2}\right)^{\mu+2n}\\&\hskip 10mm\times E_{\nu,\hskip .5mm\nu(\mu+2n)+1}(-a^{\nu}t^{\nu}),\endaligned\end{eqnarray}
where $ E_{\nu,\nu(\mu+2n)+1}(.) $ is the generalized Mittag-Leffler function defined in equation \eqref{2.3}.
\end{corollary}

By applying the results in equations \eqref{Poch} and \eqref{Poch 1}, after little simplification the  Theorems 1, 2 and 3 reduced to the following form:

\begin{corollary} If $ d>0, \nu>0, \lambda,\gamma,\mu,c,b\in\mathbb{C}, k\in\mathbb{N}$ and $\Re(\lambda)>0,\Re(\mu)>0$ then the solution of the equation
\begin{eqnarray}\aligned & N(t)-N_{0}\frac{k^{1-\mu/k-(b+1)/2k}}{\Gamma(\gamma/k)}\,{}_1\psi_2\left[
\begin{array}{cc} (\gamma/k,1);\\(\mu/k+(b+1)/k,\lambda/k),(1,1);\end{array}
t\right]\\&=-d^{\nu}{}_{0}D^{-\nu}_{t}N(t)\endaligned\end{eqnarray}

is given by the following formula

\begin{eqnarray}\aligned & N(t)=N_{0}\frac{k^{1-\mu/k-(b+1)/2k}}{\Gamma(\gamma/k)}\sum^\infty_{n=0}\frac{(-ck^{\lambda/k-1})^n}{\Gamma(\mu/k+\lambda n/k+\frac{b+1}{2k})}\frac{\Gamma(\mu+2n+1)}{(n!)^2}\\&\times\left(\frac{t}{2}\right)^{\mu+2n} E_{\nu,\hskip .5mm\mu+2n+1}(-d^{\nu}t^{\nu}),\endaligned \end{eqnarray}
where $ E_{\nu,\mu+2n+1}(.) $ is the generalized Mittag-Leffler function defined in equation \eqref{2.3}.
\end{corollary}

\begin{corollary} If $ d>0, \nu>0, \lambda,\gamma,\mu,c,b\in\mathbb{C}, k\in\mathbb{N}$ and $\Re(\lambda)>0,\Re(\mu)>0$ then the solution of the equation

\begin{eqnarray}\aligned &  N(t)=N_{0}\frac{k^{1-\mu/k-(b+1)/2k}}{\Gamma(\gamma/k)}\,{}_1\psi_2\left[
\begin{array}{cc} (\gamma/k,1);\\(\mu/k+(b+1)/k,\lambda/k),(1,1);\end{array}
d^{\nu}t^{\nu}\right]\\&-d^{\nu}{}_{0}D^{-\nu}_{t}N(t)\endaligned\end{eqnarray}

is given by the following formula

\begin{eqnarray}\aligned &  N(t)=N_{0}\frac{k^{1-\mu/k-(b+1)/2k}}{\Gamma(\gamma/k)}\sum^\infty_{n=0}\frac{(-ck^{\lambda/k-1})^n}{\Gamma(\mu/k+\lambda n/k+\frac{b+1}{2k})}\frac{\Gamma(\nu(\mu+2n)+1)}{(n!)^2}\left(\frac{d^{\nu}t^{\nu}}{2}\right)^{\mu+2n}\\&\hskip 10mm\times E_{\nu,\hskip .5mm\nu(\mu+2n)+1}(-d^{\nu}t^{\nu}),\endaligned\end{eqnarray}
where $ E_{\nu,\nu(\mu+2n)+1}(.) $ is the generalized Mittag-Leffler function defined in equation \eqref{2.3}.
\end{corollary}

\begin{corollary} If $ a>0, d>0, \nu>0 ;a\neq d; \lambda,\gamma,\mu,c,b\in\mathbb{C}, k\in\mathbb{N}$ and $\Re(\lambda)>0,\Re(\mu)>0$ then the solution of the equation

\begin{eqnarray}\aligned & N(t)=N_{0}\frac{k^{1-\mu/k-(b+1)/2k}}{\Gamma(\gamma/k)}\,{}_1\psi_2\left[
\begin{array}{cc} (\gamma/k,1);\\(\mu/k+(b+1)/k,\lambda/k),(1,1);\end{array}
d^{\nu}t^{\nu}\right]\\&-a^{\nu}{}_{0}D^{-\nu}_{t}N(t)\endaligned\end{eqnarray}

is given by the following formula

\begin{eqnarray}\aligned & N(t)=N_{0}\frac{k^{1-\mu/k-(b+1)/2k}}{\Gamma(\gamma/k)}\sum^\infty_{n=0}\frac{(-ck^{\lambda/k-1})^n}{\Gamma(\mu/k+\lambda n/k+\frac{b+1}{2k})}\frac{\Gamma(\nu(\mu+2n)+1)}{(n!)^2}\\&\hskip 10mm\times \left(\frac{d^{\nu}t^{\nu}}{2}\right)^{\mu+2n}E_{\nu,\hskip .5mm\nu(\mu+2n)+1}(-a^{\nu}t^{\nu}),\endaligned\end{eqnarray}

where $ E_{\nu,\nu(\mu+2n)+1}(.) $ is the generalized Mittag-Leffler function defined in equation \eqref{2.3}.
\end{corollary}

\section{Graphical Interpretation}

In this section we plot the graphs of main results established in equation \eqref{Th-1-2}, \eqref{Th-2-2} and \eqref{Th-3-2}.
  Graphs of the solution of the equation \eqref{Th-1-2} are depicted below for some parameter values i.e. $N_0=c=k=2;b=d=3;\mu=\nu=\gamma=1;\lambda=1,1.25,1.5,1.75,2$ in Fig. 1, Fig. 2 and Fig. 3 for time interval $t=0:1$, $t=0:2$ and $t=0:3$ respectively;
 graphs of the solution of the equation \eqref{Th-2-2} are depicted below for some parameter values i.e. $N_0=c=k=2;b=d=3;\mu=\nu=\gamma=1;\lambda=1,1.25,1.5,1.75,2$ in Fig. 4 and Fig. 5 for time interval $t=0:.05$ and $t=0:.06$ respectively.
graphs of of the solution of the equation \eqref{Th-3-2} are depicted below for some parameter values i.e. $N_0=c=k=2;b=d=3;a=\mu=\nu=\gamma=1;\lambda=1,1.25,1.5,1.75,2$ in Fig. 6 and Fig. 7 for time interval $t=0:.05$ and $t=0:.06$ respectively.
It is clear from these figures that $N(t)>0$ and the behavior of the solutions for different parameters and time interval can be studied and observed very easily. It is also observed that if we select $a=d$ in equations \eqref{Th-2-2} and \eqref{Th-3-2} give the identical solutions as we select in figures 4, 5, 6 and 7. Figures 4 and 6; 5 and 7 represents the identical solutions.

\section{Conclusion}
In this work we give a new fractional generalization of the standard kinetic equation and derived solution for the same. From the close relationship of the k-Bessel function with many special functions, we can easily construct various known and new fractional kinetic equations.

\end{document}